\documentclass[12pt]{amsart}
\usepackage{amscd,amssymb}
\usepackage[graph,frame,poly,arc]{xy}  
\usepackage{tikz}
\usepackage[plainpages,backref,urlcolor=blue]{hyperref}

\theoremstyle{plain}
\newtheorem{thm}[subsection]{Theorem}
\newtheorem{lem}[subsection]{Lemma}
\newtheorem{prop}[subsection]{Proposition}
\newtheorem{cor}[subsection]{Corollary}

\theoremstyle{definition}

\newtheorem{ex}[subsection]{Example}

\numberwithin{equation}{section}
\newcommand{\XX}{{\mathcal X}}
\newcommand{\D}{{\mathcal D}}

\newcommand{\C}{\mathbb{C}}
\newcommand{\PP}{\mathbb{P}}

\DeclareMathOperator{\rank}{rank}
\DeclareMathOperator{\brank}{brank}
\DeclareMathOperator{\crank}{crank}

\begin{document}

\title [Waring ranks of sextic binary forms via geometric invariant theory]
{Waring ranks of sextic binary forms via geometric invariant theory}

\author[Alexandru Dimca]{Alexandru Dimca$^{1}$}
\address{Universit\'e C\^ ote d'Azur, CNRS, LJAD, France and Simion Stoilow Institute of Mathematics,
P.O. Box 1-764, RO-014700 Bucharest, Romania}
\email{dimca@unice.fr}

\author[Gabriel Sticlaru]{Gabriel Sticlaru}
\address{Faculty of Mathematics and Informatics,
Ovidius University
Bd. Mamaia 124, 900527 Constanta,
Romania}
\email{gabriel.sticlaru@gmail.com}

\thanks{\vskip0\baselineskip
\vskip-\baselineskip
\noindent $^1$This work has been partially supported by the Romanian Ministry of Research and Innovation, CNCS - UEFISCDI, grant PN-III-P4-ID-PCE-2020-0029, within PNCDI III}



\begin{abstract} We determine the Waring ranks of all sextic binary forms
with complex coefficients using a Geometric Invariant Theory approach. Using the five basic invariants for sextic binary forms, our results give a rapid method to determine the Waring rank of any given sextic binary form.
In particular, we shed new light on a claim by E. B. Elliott at the end of the 19th century concerning the binary sextics with Waring rank 3.
We show that for binary forms of arbitrary degree the cactus rank, a.k.a. scheme rank, is determined by the corresponding Waring rank. Finally we determine the border ranks of all binary sextics.

\end{abstract}
 
\maketitle


\section{Introduction} 
For the general question of symmetric tensor decomposition we refer to
\cite{Guide, Car0+,Car1+, Car2+,Ci, CGLM, CS, FOS, FLOS, IK,L1,LT,O}, as well as to the extensive literature
quoted at the references in \cite{Guide} and \cite{FLOS}.
 Consider the graded polynomial ring $S=\C[x,y]$, let $S_d$ denote the vector space of homogeneous polynomials of degree $d$ in $S$, and let $f \in S_d$ be a binary form of degree $d$.  We consider the {\it Waring decomposition}
\begin{equation}
\label{eq1}
(\D) \ \  \  \ \  f=\ell_1^d + \cdots +\ell_r^d,
\end{equation}
where $\ell_j \in S_1$ are linear forms in $x,y$, and $r$ is minimal, in other words $r=\rank f$ is the {\it Waring rank}  of $f$. The {\it border rank} $\brank f$ of the form $f$ is the smallest integer $r'$ such that $f$ is in the closure of the set of binary forms of Waring rank $r'$ in $S_d$.
 Let $Q=\C[X,Y]$,  where $X=\partial_x$ and $Y=\partial_y$. Then $Q$ is the ring of differential operators with constant coefficients and acts on $S$ in the obvious way. For a binary form $f \in S$, we consider the ideal
of differential operators in $Q$ killing $f$, namely
\begin{equation}
\label{eq2}
Ann(f)=\{q \in Q : q \cdot f=0\},
\end{equation}
also denoted by $f^{\perp}$ and called the {\it apolar ideal} of $f$.
The {\it cactus rank} $\crank f$ of the form $f$ is the minimal length of a $0$-dimensional subscheme $\Gamma$ in $\PP^1=\PP(Q_1)$, defined by a homogeneous ideal $I_{\Gamma}$ in $Q$ satisfying $I_{\Gamma} \subset Ann(f)$, see \cite{BBM} for the general setting. Note that sometimes, as in \cite{IK} for instance, the cactus rank is called {\it scheme rank}.

It is clear that these three ranks of a non zero binary form $f$ of degree $d$ depend only on the corresponding class $[f]$ in the projective space $\PP(S_d)$, and even on the corresponding $SL_2(\C)$-orbit of $[f]$ in $\PP(S_d)$. In the sequel, to simplify the notation, we denote $[f]$ simply by $f$ and talk of the $SL_2(\C)$-orbit of $f$ in $\PP(S_d)$.

The Waring ranks of binary forms of degree $\leq 5$ are rather well understood, see \cite{BM,DStToh}.
In this note we discuss the Waring ranks of binary sextics using invariant theory of such forms, which is recalled briefly in section 3. There are five basic invariants for sextic binary forms
\begin{equation}
\label{eqI0}
f=a_0x^6+6a_1x^5y+15a_2x^4y^2+20a_3x^3y^3+15a_4x^2y^4+6a_5xy^5+a_6y^6,
\end{equation}
classically denoted by $I_2(f),I_4(f),I_6(f),I_{10}(f)$ and $I_{15}(f)$. 
In classical terminology, $I_2(f)$ is the {\it sixth transvectant of
$f$ over itself}, transvection being one of the main operations in Gordan's proof of
finiteness for invariants of binary forms. Note that $I_4(f)$ is the determinant of the {\it catalecticant matrix} $C(f)_3$ of $f$ from Example \ref{ex1}.

It was claimed by E. B. Elliott at the end of the 19th century that $I_4(f)$ vanishes precisely when $f$ can be expressed as a sum of three sixth powers, see \cite[Exercise 10] {El}, pp. 276. This fact is also quoted in modern literature, see for instance
 \cite{H}, the discussion before Proposition 2.1. We clarify this claim in Theorem \ref{thm2}, by showing that the Waring rank of a sextic $f$ satisfying $I_4(f)=0$ can be as high as 5.

Our results can be summarized as follows. 
From the invariant theory view-point, the sextic binary forms can be divided into 3 classes as follows, see Proposition \ref{prop1} and Theorem \ref{thm1}.

\begin{enumerate}

\item The forms in the Null Cone, that is the forms on which all the above invariants vanish. These are exactly the sextic forms having a root of multiplicity at least 4. For them the Waring ranks are easily determined, see Lemmas \ref{lemM6}, \ref{lemM5} and \ref{lemM4}.

\item The forms which are semi-stable, but not stable. These are exactly the sextic forms having a root of multiplicity 3 and they have all Waring rank 4, see Proposition \ref{propM3}.

\item The stable forms, which are exactly the sextic forms having at most double roots. The $SL_2(\C)$-orbit of such a form is determined by the point
$$(I_2(f):I_4(f):I_6(f):I_{10}(f)) \in \XX=\PP(1,2,3,5) \setminus \{(10:1:8^3:0)\},$$
where $\PP(1,2,3,5)$ denotes the weighted projective space, with coordinates $z_1,z_2,z_3,z_5$ and weights $wt(z_j)=j$, for $j=1,2,3,5$.
The open set in $\XX$ given by $z_2 \ne 0$ consists of forms of Waring rank 4, and this is indeed the generic situation, see Proposition \ref{propS1}. 

In the hyperplane $z_2=0$, which can be identified with the weighted projective plane $\PP(1,3,5)$, the complement of the rational curve
$$C: z_1^3+z_3=z_2=0$$
consists of forms of rank 3. There is a special point on this curve,  
$$p=(1:0:-1:1) \in C,$$
such that the corresponding forms are exactly those of Waring rank 2.
Finally, the pointed curve $C \setminus \{p\}$ correspond to $SL_2(\C)$-orbits of forms with Waring rank 5. All these claims are proved in Theorem \ref{thm2}.
\end{enumerate}

An interesting fact is that the presence of a root of multiplicity $\geq 3$ determines the Waring rank of $f$. The following result is essentially well known, but we include a proof for the reader's convenience, collecting facts from Lemmas \ref{lemM6}, \ref{lemM5}, \ref{lemM4} and Proposition \ref{propM3}. A similar property holds for the border ranks of sextic binary forms, see Corollary \ref{corGITborder}.

\begin{cor}
\label{corGIT}
Let $f$ be a non zero binary sextic having a  root of highest multiplicity $m \geq 3$. Then the following hold.

\begin{enumerate}

\item  $m=6$ if and only if  $\rank f=1$.

\item $m=5$ if and only if  $\rank f=6$.

\item If $m=4$, then $\rank f=5$.

\item If $m=3$, then $\rank f=4$.

\end{enumerate}

\end{cor}

On the other hand, when the maximal multiplicity of a root of $f$ is at most 2, then there are a lot of possibilities, which we list now in our first main result, and prove 
after  Example \ref{exGIT}. To separate all these cases, the invariant theory approach plays a key role in our opinion.
\begin{thm}
\label{thmcorGIT}
Let $f$ be a binary sextic having at most double roots. The following cases may occur.

\begin{enumerate}

\item If $f$ has three double roots, then $\rank f=4$.

\item If $f$ has two double roots, then $3 \leq \rank f \leq 4$, and both cases are possible.

\item If $f$ has a single double root, then $3 \leq \rank f \leq 5$, and the three cases are possible.

\item If $f$ has distinct roots, then $2 \leq \rank f \leq 5$, and the four cases are possible.

\end{enumerate}

\end{thm}
Note that, given a sextic binary form $f$ as in Theorem \ref{thmcorGIT}, we can 
decide rapidly using our results the Waring rank  of a binary sextic $f$ by just computing the corresponding vector 
$$((I_2(f),I_4(f),I_6(f),I_{10}(f)).$$
Our results allow us to determine the border rank of a binary sextic $f$, see Theorem \ref{thmBRK}. On the other hand, in the case of binary forms of any degree, Sylvester's Theorem implies that the Waring rank determines the cactus rank by a simple formula, see Corollary \ref{corCRK}. Therefore, our results determine the three invariants $\rank f$, $\brank f$ and $\crank f$ of any  binary sextic $f$.

\medskip

 Computations with CoCoa \cite{Co} and Singular \cite{Sing} played a key role in our results.
We would like to thank Nicolae Manolache for useful discussions concerning the invariant theory of sextic binary forms and Bruce Reznick for his help in proving the claim (3) in Theorem \ref{thmBRK}. We would also like to thank the referee for the very careful reading of our manuscript and for his suggestions to improve the presentation and to add the analysis of the border rank and the cactus rank of binary sextics.

\section{Sylvester's Theorem} 

For the reader's convenience, we recall in this section Sylvester's Theorem which plays a key role in our discussion. Note that the apolar ideal $Ann(f)$ defined in \eqref{eq2} is a graded ideal, whose degree $s$ homogeneous component is given by $\ker [f]$, where
\begin{equation}
\label{eq3}
[f]: Q_s \to S_{d-s}
\end{equation}
is the morphism $g \mapsto g \cdot f$. The matrix of this linear map with respect to the obvious monomial bases of the vector spaces $Q_s$ and $ S_{d-s}$ is called the {\it catalecticant matrix} $C(f)_s$ of $f$ in degree $s$.

\begin{ex}
\label{ex1} 
As an example, if we take $d=6$ and $s=3$, then for $f$ as in \eqref{eqI0} we get
\begin{center}
$$C(f)_3= 2^{12}3^{6}5^{4} \left(
  \begin{array}{ccccccc}
     a_0 & a_1& a_2 & a_3 \\
     a_1 &a_2& a_3 & a_4\\
    a_2& a_3& a_4 & a_5 \\
    a_3& a_4& a_5 & a_6 \\
   \end{array}
\right).$$
\end{center}

\end{ex}

One has the following well known result, see for instance\cite[Lemma 2.2]{DStToh}.
\begin{lem}
\label{lemS}
\begin{enumerate}
The graded ideal $Ann(f) \subset Q$ of the binary form $f$ of degree $d$ satisfies the following.

\item $Ann(f)_0 \ne 0$ if and only if $f=0$.

\item $Ann(f)_0 = 0$ and $Ann(f)_1 \ne 0$ if and only if $f=x^d$  after a linear change of coordinates.

\item $Ann(f)_1= 0$ and $Ann(f)_2= \C \ell^2$ for some $\ell \in Q_1$
if and only if $f=x^{d-1}y$  after a linear change of coordinates.

\end{enumerate}

\end{lem}

The following result goes back to Sylvester \cite{Sy}. See also \cite{CS}.

\begin{thm}
\label{thmS}
For a binary form $f$ of degree $d$,
the apolar ideal $Ann(f)$ is a complete intersection, namely there are two binary forms $g_1$ and $g_2$ in $Q$ such that $Ann(f)=(g_1,g_2)$.
The degrees $d_j$ of $g_j$ for $j=1,2$ satisfy $d_1+d_2=d+2$.
Moreover, if we assume $d_1 \leq d_2$, then the Waring rank $r=\rank f$ is determined as follows.
\begin{enumerate}

\item If the binary form $g_1$ has no multiple factors, then $r=d_1$.

\item Otherwise, $r=d_2$.

\end{enumerate}

\end{thm}
According to Lemma \ref{lemS}, the interesting case is 
$$2 \leq d_1 \leq d_2.$$
In this case we have the following, see also the Introduction in \cite{BT}.
\begin{thm}
\label{thmAH}
If the binary form $f$ of degree $d$ satisfies $\rank f\geq 2$, then
$$\rank f \leq d$$
and the equality holds if and only if $f=x^{d-1}y$  after a linear change of coordinates.
Moreover, for a generic binary form $f$ of degree $d$ one has
$$\rank f = \left \lfloor \frac{ d+2}{2} \right \rfloor=  \left \lceil \frac{ d+1}{2} \right \rceil.$$
\end{thm}
\proof
The first claim follows from Lemma \ref{lemS}(3) and
Sylvester's Theorem \ref{thmS}. The second claim is a special case of 
Alexander-Hirschowitz results in \cite{AH}.
\endproof
In particular, we see that for a generic sextic binary form $f$ one has
$$g= \rank f=4.$$
Generally speaking, the forms in $n$ variables of a given degree $d$, having a {\it Waring rank strictly higher than the Waring rank of a generic form} are very mysterious, and they are studied in full generality (any degree $d$ and any number of variables $n$) in \cite{B2}. When $n=2$, then the sets of binary forms of a given Waring rank can be described in terms of joins of algebraic subvarieties in $\PP(S_d)$, see \cite[Section (4.4)]{B2}.
When $n=2$ and $d=6$, such high rank forms can have either
$\rank f=6$, and this case is very easy to understand, see Lemma \ref{lemM5}, or $\rank f=5$, and here there are several subtle possibilities, see 
Lemma \ref{lemM4} and Theorem \ref{thm2}, Case (2) (b).

As a consequence of Sylvester's Theorem \ref{thmS} we have the following result, showing that the Waring rank $\rank f$ determines the cactus rank $\crank f$ for binary forms of any degree.

\begin{cor}
\label{corCRK}
Let $f$ be a non zero binary form of degree $d$. With the notation from Theorem \ref{thmS}, one has the following relations between the cactus rank and the Waring rank of the binary form $f$.
\begin{enumerate}

\item  If $r=\rank f \leq  \frac{ d+2}{2} $, then $\crank f= \rank f$.

\item If $r=\rank f >  \frac{ d+2}{2} $, then $\crank f= d+2 -\rank f$.

\end{enumerate}

\end{cor}
\proof
A $0$-dimensional subscheme $\Gamma$ in $\PP^1=\PP(Q_1)$ is defined by a principal  homogeneous ideal $I_{\Gamma}=(g)$ in $Q$ and its length is just $\deg g$. 
With the notation from Theorem \ref{thmS}, it is clear that the scheme $\Gamma_1$ defined by the principal ideal $(g_1) \subset Q$ is a minimal length $0$-dimensional subscheme satisfying $I_{\Gamma_1} \subset Ann(f)$. The case (1) in Corollary \ref{corCRK} corresponds to the case (1) in Theorem \ref{thmS}, and the case (2) in Corollary \ref{corCRK} corresponds to the case (2) in Theorem \ref{thmS}. 
\endproof

\section{Invariant Theory of sextic binary forms}

Following \cite[Section 10.2]{Do}, let
$Hyp_6(1)=\PP(S_6)=\PP^6$
be the projective space parametrizing all nonzero binary forms $f$ of degree $d=6$, written as in \ref{eqI0}. 

Consider the natural $SL_2(\C)$-action on $Hyp_6(1)$ induced by the coordinate changes of $\PP^1$. Then the subset of stable points $Hyp_6(1)^s$ (resp. the subset of semi-stable points $Hyp_6(1)^{ss}$) in $Hyp_6(1)$ with respect to this $SL_2(\C)$-action
action is given by the forms having no multiple roots of multiplicity
$\ge 3$ (resp. $>3$), see \cite[Theorem 10.2]{Do}.
The categorical quotient
$$\pi : Hyp_6(1)^{ss} \to C_6(1)=Hyp_6(1)^{ss}//SL_2(\C)$$
can be described explicitly since we know the algebra of $SL_2(\C)$-invariant polynomials $A=Pol(S_6)^{SL_2(\C)}$ defined on the vector space $S_6$. One has in fact
$$C_6(1)=Projm(A),$$
the projective variety associated to the graded finitely generated $\C$-algebra $A$ in the usual way,
and the mapping $\pi$ is given by a set of generators of the $\C$-algebra $A$, see \cite[Section 10.2]{Do}. 
It is known that the $\C$-algebra $A$ is generated by five invariants, classically denoted by $I_2,I_4,I_6,I_{10}$ and $I_{15}$, satisfying a unique relation
$$I_{15}^2=F(I_2,I_4,I_6,I_{10})$$
for some polynomial $F$. However, the subalgebra $A^{(2)}$ of $A$ generated by invariants of even degree is a polynomial algebra and using this, as noticed in \cite[Section 10.2]{Do}, we can identify $Projm(A)$ with the weighted projective space 
$$Projm(A^{(2)}) =\PP(2,4,6,10)=\PP(1,2,3,5).$$ 
Then the canonical projection is given by
\begin{equation}
\label{eqI1}
\pi: Hyp_6(1)^{ss} \to   \PP(1,2,3,5),  \   \     \pi(f)=(I_2(f):I_4(f):I_6(f):I_{10}(f)).
\end{equation}
The first two basic invariants can be chosen  as follows, see \cite{Bo,Cl,Ig,H,T},
\begin{equation}
\label{eqI2}
I_2(f) = -10a_3^2 + 15a_2a_4 - 6a_1a_5 + a_0a_6,
\end{equation}
and
\begin{equation}
\label{eqI4}
I_4(f) = \det \left(
  \begin{array}{ccccccc}
     a_0 & a_1& a_2 & a_3 \\
     a_1 &a_2& a_3 & a_4\\
    a_2& a_3& a_4 & a_5 \\
    a_3& a_4& a_5 & a_6 \\
   \end{array}
\right)=
\end{equation}
$$=a_3^4-3a_2a_3^2a_4 + a_2^2a_4^2+ 2a_1a_3a_4^2-a_0a_4^3+ 2a_2^2a_3a_5 
-2a_1a_3^2a_5 -2a_1a_2a_4a_5 + $$
$$+2a_0a_3a_4a_5
+a_1^2a_5^2-a_0a_2a_5^2-a_2^3a_6 + 2a_1a_2a_3a_6 -a_0a_3^2a_6 -a_1^2a_4a_6 + a_0a_2^4a_6.$$
To define $I_6(f)$, we proceed as follows. First we set 
$$ b_0=6(a_0a_4-4a_1a_3+3a_2^2), \ b_1 = 3(a_0a_5-3a_1a_4 + 2a_2a_3),$$
$$ \ b_2 = a_0a_6 -9a_2a_4 + 8a_3^2, \ b_3 = 3(a_1a_6 -3a_2a_5 + 2a_3a_4), \ b_4 = 6(a_2a_6- 4a_4a_5 + 3a_4^2).$$

Then we define
\begin{equation}
\label{eqI5}
I_6(f)=\det \left(
  \begin{array}{ccccccc}
     b_0 & b_1& b_2  \\
     b_1 &b_2& b_3 \\
    b_2& b_3& b_4  \\
   \end{array}
\right).
\end{equation}
Finally we set
$$I_{10}(f)=6^{-10} \Delta(f),$$
where $\Delta(f)$ is the discriminant of the sextic form $f$ in \ref{eqI0}.
One has the following result. \begin{prop}
\label{prop1}
For a binary sextic $f \in Hyp_6(1)$ we have the following.
\begin{enumerate}

\item $I_{10}(f)=0$ if and only if $f$ has a multiple root.

\item $\pi(f)=(10:1:8^3:0)$ if $f$ has a triple root.

\item  $I_{2}(f)=I_{4}(f)=I_{6}(f)=I_{10}(f)=0$ if and only if  $f$ has a root of multiplicity at least 4.
\end{enumerate}
\end{prop}

\proof The first claim is clear. For the second claim, assume that the root of multiplicity 3 is $x=0$. Then $a_4=a_5=a_6=0$ and $a_3 \ne 0$.
The claim follows using the formulas for the invariants and the equality
$$(-10:1:-8^3:0)=(10:1:8^3:0)$$
in $\PP(1,2,3,5)$.
For the last claim we refer to see  \cite[Lemma 3.3]{KSV}.
\endproof

We summarized now  the results we need in the sequel.
\begin{thm}
\label{thm1}
\begin{enumerate}

\item $Hyp_6(1) \setminus Hyp_6(1)^{ss}$ coincides with the Null Cone $N_6$ given by 
$$N_6=\{f \in Hyp_6(1): \ I_{2}(f)=I_{4}(f)=I_{6}(f)=I_{10}(f)=0\}.$$
Moreover the binary sextic $f\in N_6$ if and only if $f$ has a root of multiplicity at least 4.

\item The set of semi-stable, but not stable sextics is given by
 the sextics $f$ which have a root of multiplicity exactly 3. All these forms are mapped to the point
$(10:1:8^3:0) \in \PP(1,2,3,5)$ under the quotient map $\pi$.

\item $\pi: Hyp_6(1)^{s} \to \XX=\PP(1,2,3,5) \setminus \{(10:1:8^3:0)\}$ is a geometric quotient for binary sextics with zeroes of multiplicity $\leq2$.
In particular,  if the sextic $f$ has at most roots of multiplicity 2, then $f\in Hyp_6(1)^{s}$  and hence the $SL_2(\C)$-orbit
$SL_2(\C)\cdot f$ coincides with $\pi^{-1}(\pi(f)).$

\end{enumerate}

\end{thm}

\proof
The  claim (1) follows from Proposition \ref{prop1} (3), in view of the description of semi-stability given above. The claim (2) follows from 
Proposition \ref{prop1} (2) and the description of stability given above.
The claim (3) is a general property in geometric invariant theory, see
for instance \cite[Theorem 8.1]{Do} or \cite[Proposition 2.3 (2)]{H}, but notice the difference of the choice for the basic invariants in \cite{H}.
\endproof

In particular, Theorem \ref{thm1} (3) implies the following classical result, see \cite{Cl,Bo,Ig} as well as
 \cite[Corollary 3.13]{KSV}.
\begin{cor}[Clebsch-Bolza-Igusa]
\label{cor1}
Two binary sextics $f,g \in \PP(S_6)$ without multiple roots are in the same $SL_2(\C)$-orbit if and only if
$$\pi(f)=\pi(g)$$
in $\PP(1,2,3,5)$.
\end{cor}

\section{The Waring ranks of sextic binary forms}

We determine in this section the Waring rank of a non-zero sextic binary form $f \in Hyp_6(1)$ using mainly the invariant theory recalled in Section 3. First we study the forms in the Null Cone $N_6$, where invariant theory cannot help. The first four results below are essentially well known, we state them to have a clear reference, and include a proof for the reader's convenience.

\begin{lem}
\label{lemM6}  The sextic  $f$ has a root of multiplicity 6 if and only if $\rank f=1$.
\end{lem}
\proof
Obvious from the definition of the Waring rank.
\endproof

\begin{lem}
\label{lemM5}
  The sextic $f$ has a root of multiplicity 5  if and only if $\rank f=6$.
\end{lem}
\proof  
A binary form $f \in Hyp_6(1)$ has a root of multiplicity 5 if and only if  
$f \in SL_2(\C)\cdot x^5y$. As recalled in Lemma \ref{lemS}, this happens exactly when 
$$\rank C(f)_1 =2 \text{ and } \rank C(f)_2 \leq 2,$$
and this is the case precisely when the Waring rank of $f$ is 6 by Theorem \ref{thmAH}.
\endproof

\begin{lem}
\label{lemM4}
 If $f$ has a root of multiplicity 4, then $\rank f=5$.
\end{lem}
\proof   
  If a binary form $f \in Hyp_6(1)$ has a root of multiplicity 4,  then it is known that $\rank f \geq 5$, see \cite[Theorem 3.1]{Tok}. On the other hand, we know by Theorem \ref{thmAH} that $\rank f \leq 6$, and by the previous case we know that the equality $\rank f = 6$ cannot hold when $f$ has a root of multiplicity 4. Algebraically, this case corresponds to the conditions
$$\rank C(f)_2 =3  \text{ and } I_{2}(f)=I_{4}(f)=I_{6}(f)=I_{10}(f)=0.$$
\endproof

This ends the discussion of the forms in the Null Cone. Next we treat the semi-stable forms which are not stable.

\begin{prop}
\label{propM3}
 If $f$ has  a root of multiplicity 3, then  $\rank f=4$.
 \end{prop}
 \proof
In this case, $f$ is in the $SL_2(\C)$-orbit determined by one of the following 3 forms: $f_1=x^3y^3$, $f_2=x^3y^2(x+y)$ or $f_3=x^3y(x+y)(x+ty)$, for some $t \ne 0,1$ according to the multiplicities of the other roots of $f$.
If  $f \in SL_2(\C)\cdot f_j$, then $\rank f= \rank f_j$ and hence we have to determine the Waring rank of $f_j$, for $j=1,2,3$. 
Note that according to Lemma \ref{lemS}, $Ann(f_j)_1=Ann(f_j)_2=0$.
Then a direct computation using Example \ref{ex1} shows that 
$$\det C(f_j)_3 \ne 0$$
and hence $Ann(f_j)_3=0$. It follows that in Theorem \ref{thmS} applied to $f_j$, we have
$$4 \leq d_1 \leq d_2 \text{ and } d_1+d_2=8.$$
Therefore 
$$\rank f_j=d_1=d_2=4$$
is the only possibility, for $j=1,2,3$.
\endproof

From now on we consider only stable binary sextics, i.e. sextics with at most double roots, and determine their Waring ranks using the fact their orbits $SL_2(\C)\cdot f$ are determined by $\pi(f)$.

\begin{prop}
\label{propS1}
If $f$ is stable, then  $I_4(f) \ne 0$ if and only if  $\rank f=4$.
 \end{prop}
 \proof
This is the generic case, compare to the last claim in Theorem \ref{thmAH}. If
$$I_4(f)=\det C(f)_3 \ne 0,$$
it follows that $Ann(f)_k=0$ for all $k \leq 3$. Then we conclude that $\rank f =4$ as in the proof of Proposition \ref{propM3} using Theorem \ref{thmS}. Conversely, assume now that $f$ is stable and it has Waring rank 4. Then $f$ is $SL_2(\C)$-equivalent to a normal form
$$f(u,v,t)=ux^6+vy^6+(x+y)^6+(x+ty)^6,$$
where $u, v \in \C^*$ and $t \in \C \setminus \{0,1\}$.
A direct computation shows that
$$I_4(f(u,v,t))=uv(t-1)^2t^2 \ne 0.$$
Hence $I_4(f)$, which is equal to a non-zero multiple of $I_4(f(u,v,t))$, is also non-zero.
\endproof

\begin{thm}
\label{thm2}
 If $f$ is stable and $I_4(f) = 0$, then one and only one of the following cases may occur. 
 
\begin{enumerate}

\item   $I_2(f)^3+I_6(f) \ne 0$ and then $\rank f=3$, 

\item $I_2(f)^3+I_6(f) = 0$ and then

\begin{enumerate}

\item either $\rank f=2$ and $\pi(f)= (1:0:-1:1)$, or

\item $\rank f=5$,  and $\pi(f)= (0:0:0:1)$ or $\pi(f)=(1:0:-1:t)$ for some $t\in \C$, $t \ne 1$.

\end{enumerate}

\end{enumerate}

 \end{thm}

 \proof
Since
$I_4(f)=\det C(f)_3= 0,$
it follows that $Ann(f)_3\ne 0$. Moreover, since $f$ is stable, we know that $Ann(f)_1=0$.

Suppose first that $Ann(f)_2 \ne 0$. Then, using stability and Lemma \ref{lemS}(3), we deduce from Theorem \ref{thmS} that $\rank f =2$.
In this case $f$ is $SL_2(\C)$-equivalent to $f_1=x^6+y^6$ and hence
$$\pi(f)=\pi(f_1)= (1:0:-1:1).$$
This gives the subcase (a) in our Theorem.

Assume now that $Ann(f)_2 =0$, and hence $\rank f \geq 3$.
We know that $Ann(f)_3\ne 0$ and hence $d_1=\deg(g_1)=3$. In fact $\dim Ann(f)_3 =1$, since otherwise in Theorem \ref{thmS} we would have $d_1=d_2=3$ and hence $d_1+d_2 <d+2=8$, a contradiction.  Let $\delta$ be a generator of $Ann(f)_3$.
Up to a linear change in the varables $X,Y$, there are 3 cases to consider.
\medskip

\noindent {\bf Case A. $\delta=X^3$} \\
Then $Xf=0$ implies that $f$ contains only monomials $x^ky^{6-k}$ with
$0 \leq k \leq 2$. This implies that $f$ is divisible by $y^4$, in contradiction with the stability of $f$. Hence this case cannot occur.

\medskip

\noindent {\bf Case B. $\delta=X^2Y$} \\
Then $Xf=0$ implies that the coefficients of $f$ given in \eqref{eqI0} satisfy
$$a_1=a_2=a_3=a_4=0.$$
Since $f$ is stable, $a_0 \ne 0$, otherwise $y=0$ would be a root of multiplicity $\geq 5$. Note also that $a_5 \ne 0$, since otherwise $\rank f=2$ and hence $Ann(f)_2 \ne 0$, a contradiction.
It follows that
$$I_2(f)=a_0a_6, \ I_4(f)=0, \ I_6(f)=-(a_0a_6)^3 \text{ and } I_{10}(f)= a_0^4(a_0a_6^5-5^5a_5^6),$$
and hence
$$\pi(f)=(a_6:0:-a_6^3: a_6^5-5^5a_5^6/a_0).$$
This gives $\pi(f)=(0:0:0:1)$ if $a_6=0$ and
$$\pi(f)=(1:0:-1:1- \frac{5^5a_5^6}{a_0a_6^5}) \ne (1:0:-1:1),$$
if $a_6 \ne 0$. This yields the subcase (b).

\medskip

\noindent {\bf Case C. $\delta=XY(X+Y)$} \\
Then Theorem \ref{thmS} implies that $\rank f =3$. Any sextic of Waring rank 3 is $SL_2(\C)$-equivalent to the normal form
$$f(u,v)=ux^6+vy^6+(x+y)^6$$
for some $u,v \in \C^*$. A direct computation shows that
$$I_2(f(u,v))=uv+u+v$$
and
$$I_6(f(u,v))=-u^3v^3-3u^3v^2-3u^3v-u^3-3u^2v^3+48u^2v^2-3u^2v-3uv^3-3uv^2-v3.$$
It follows that
$$I_6(f(u,v)) + I_2(f(u,v))^3 = 54u^2v^2 \ne 0.$$
This gives the case (1) in our Theorem.
\endproof
It was shown in \cite{BM} that the Waring rank of a binary form which is a binomial is easy to determine just looking at the corresponding exponents. The following example shows that for the trinomial binary forms, the determination of the Waring rank is much more complicated.
\begin{ex}
\label{exGIT}
Consider the family of binary sextics
$$f_t=x^2y^2(x+y)(x+ty)=x^4y^2+(1+t)x^3y^3+tx^2y^4,$$
where $t \ne 0$, having at least two double roots. The involution $(x,y) \mapsto (y,x)$ shows that
$f_t$ and $f_{\frac{1}{t}}$ have the same Waring rank.

\medskip

\noindent{\bf Case 1: $t \ne 1$.} Then $f_t$ has 2 double roots and 2 simple roots. Moreover, for any binary sextic $f$ which has 2 double roots and 2 simple roots
there is a value of the parameter $t$ such that $f$ and $f_t$ are $SL_2(\C)$-equivalent.
To determine the Waring rank of $f_t$ following the method of this paper, we compute the invariants for $f_t$. Note that, by replacing $f_t$ by $60f_t$, we can take
$$a_0=a_1=a_5=a_6=0, \ a_2 =a_4=4  \text{ and } a_3 =3(1+t).$$ 
It follows that
$$I_2(f_t)=-30(3t^2-2t+3), \ 
I_4(f_t)= 81t^4-108t^3-122t^2-108t+81, $$
and
$$
I_6(f_t)=-373248(t^4+2t^3-2t^2+2t+1)(t^2-4t+1).$$
The equation $I_4(f_t)=0$ has two real roots $t_1$ and $t_2=t_1^{-1}$, where
$$t_1= \frac{3+4 \sqrt 5-2 \sqrt{6\sqrt 5+2}}{9},$$
and two complex roots $t_3$ and $t_4=t_3^{-1}$, where
$$t_3= \frac{3-4 \sqrt 5-2 i\sqrt{6\sqrt 5-2}}{9}.$$
One can easily check that these roots are not roots for the polynomial
$$P(f_t)=I_2(f_t)^3+I_6(f_t).$$
Using Proposition \ref{propS1} and Theorem \ref{thm2} we conclude that
\begin{enumerate}

\item either $t \notin \{0,1,t_1,t_2,t_3,t_4\}$ and then $\rank f_t= 4$, or

\item $t \in \{t_1,t_2,t_3,t_4\}$ and then $\rank f_t= 3$.

\end{enumerate}

\medskip

\noindent{\bf Case 2: $t =1$.} Then $f_t$ has 3 double roots. Moreover, for any binary sextic $f$ which has 3 double roots,
 $f$ and $f_1$ are $SL_2(\C)$-equivalent.
The above computations show that $I_4(f_1) \ne 0$, and hence  Proposition \ref{propS1} implies that $\rank f_1=4$.
\end{ex}

\medskip

{\bf Proof of Theorem \ref{thmcorGIT}.}
The claims (1) and (2) follows from Example \ref{exGIT}.
For the claim (3), note that $3 \leq \rank f \leq 5$ follows from 
Proposition \ref{propS1} and Theorem \ref{thm2}. Indeed, the forms of rank 2 have clearly distinct roots. We have to check that all these 3 values for $\rank f$ may really occur.
For the case of rank 3, one can take
$$f=(x+y)^6+(x-y)^6-2x^6=y^2(30x^4+30x^2y^2+2y^4).$$
For the case of rank 4, one can take  
$$f=(x+y)^6+(x-y)^6-2x^6+28y^6=30y^2(x^4+x^2y^2+y^4),$$
where $a_0=a_1=a_3=a_5=0$, $a_2=a_4=2$, $a_6=30$ and hence
$I_4(f) \ne 0$.
For the rank 5, the proof of Case B in Theorem \ref{thm2} shows that we can take
$$f= x^6+6xy^5+5y^6=(x+y)^2(x^4-2x^3y+3x^2y^2-4xy^3+5y^4).$$
The last claim is clear, since the condition of distinct roots is just $I_{10}(f) \ne0$. For instance, to have an example of rank 5, we can take
$$f= x^6+6xy^5+y^6,$$
as implied by the proof of Case B in Theorem \ref{thm2}.
\endproof

\section{The border ranks of binary sextics} 

In this section we determine the border ranks of all binary sextics.

\begin{thm}
\label{thmBRK}
 Let $f$ be a non zero binary sextic form. Then its border rank $\brank f$ is determined as follows.
 
\begin{enumerate}

\item   The sextic $f$ has a root of multiplicity $6$ if and only if $\brank f=1$.

\item If the sextic $f$ has a root of multiplicity $5$, then $ \brank f =2$.

\item If the sextic $f$ has a root of multiplicity $4$,  then $\brank f=3$.

\item If the sextic $f$ has a root of multiplicity $3$, then $\brank f =4$.

\item If the sextic $f$ is stable with $I_4(f) \ne 0$, then $\brank f =4$.

\item If the sextic $f$ is stable with $I_4(f) =0$ and $I_2^3(f)+I_6(f) \ne 0$, then $\brank f =3$.

\item If the sextic $f$ is stable with $I_4(f) =0$ and $I_2^3(f)+I_6(f) = 0$, then one of the following cases can occur.

\begin{enumerate}

\item either $I_2^5(f)-I_{10}(f) = 0$ and then $\brank f=2$, or

\item $I_2^5(f)-I_{10}(f) \ne 0$ and then $\brank f =3$.

\end{enumerate}

\end{enumerate}

 \end{thm}
 
 \proof
The forms $f$ having a root of multiplicity $6$ form a single, closed, 1-dimensional  $SL_2(\C)$-orbit  in $\PP(S_6)$, and hence the claim (1) is clear. 

To prove the claim (2), we use the family of sextics 
$$f_t=(x+ty)^6-x^6$$
with $t \in \C^*$ small. It is clear that $\rank f_t=2$ and that the class of $f_t$ in $\PP(S_6)$ coincides with the class of
$$f'_t=y(6x^5+15tx^4y+20t^2x^3y^2+15t^3x^2y^3+6t^4xy^4+t^5y^5).$$
When $t$ tends to $0$, this class has as limit the class of $x^5y$, and this proves our claim.

The generic binary sextic has rank 4 as follows from Theorem \ref{thmAH}, and hence for any binary sextic $f$ one has $\brank f \leq 4$.
On the other hand, the binary sextics $f$ with $2 \leq \rank f \leq 3$ are of two types, namely exactly those described in our Theorem \ref{thm2}. Hence to decide the border rank of a binary sextic $f$ comes down to deciding if in small neighborhoods of $f$ one can find forms of these two types.
Since for them the invariant $I_4$ vanishes, it follows that any binary sextic $f$ with $I_4(f) \ne 0$ has $\brank f=4$. In particular, this proves our claim
(4) in view of Theorem \ref{thm1}, (2) and (5).

The claim (6) is also clear by Theorem \ref{thm2} (1) and the description of sextics with Waring rank $\leq 2$.

The claim (7) (a) is obvious, since such a binary sextic $f$ has $\rank f =2$. Such sextics form a single, 3-dimensional  $SL_2(\C)$-orbit  in $\PP(S_6)$. 
To prove the claim (7) (b), note first that the set of forms which occur in Theorem \ref{thm2} (2) union with the Null Cone $N_6$ is defined by the equations $I_4(f)=I_2(f)^3+I_6(f)=0$.  Theorem \ref{thm2} (2) and the description of the Null Cone given in Theorem \ref{thm1} (1) imply that the dimension of this union is 4, and hence the hypersurfaces $W:I_4(f)=0$ and $W': I_2(f)^3+I_6(f)=0$ in $\PP(S_6)$ have no common irreducible component. It follows that for any irreducible component $W_0$ of $W$, the intersection $W_0'=W_0 \cap W'$ is a proper closed subset in $W_0$. Hence $W_0'$ is in the closure of the complement $W_0 \setminus W_0'$. This complement is formed of sextics with Waring rank 3 by Theorem \ref{thm2} (2), while all the forms considered in the claim 7 (b) are contained in subsets of type $W_0'$.

For the claim (3), we assume first that $f$ has a root of multiplicity 4 and two simple roots. Then we have  $\brank f \geq 3$. Indeed, the
binary forms with $\rank f =1$ form a 1-dimensional $SL_2(\C)$-orbit  in $\PP(S_6)$ and those with $\rank f =2$ form a 3-dimensional $SL_2(\C)$-orbit  in $\PP(S_6)$. Hence the 3-dimensional orbit formed by the binary forms with a root of multiplicity 4 and 2 simple roots cannot be in the closure of either. Indeed, the closure of an orbit contains only strictly lower dimensional orbits. To prove the converse inequality $\brank f \leq 3$, we consider the form $f_0=x^4(x^2-15y^2)$ and the small deformation $f_t$ of it given by the formula \eqref{eqI0} with the choices
$a_0=1, a_2=-1$, $a_1=a_3=a_5=0$, $a_4=t$ and $a_6=-t^2$, with $t \in \C^*$ small.
A direct computation shows that
$$
I_4(f_t)=0 \text{ and } 
I_2(f_t)^3+I_6(f_t) = -216t^5+432t^4-216t^3.$$
For $t \ne 0$ but small enough, it follows that $\rank f_t=3$. Indeed, Theorem \ref{thm1} implies that $f_t$ is stable and  then Theorem \ref{thm2} (1) can be used to conclude that $\rank f_t=3$.

Consider now the case when the form  $f$ has a root of multiplicity 4 and one double root. Then  the inequality $\brank f \leq 3$ follows from the previous discussion and the fact that the 2-dimensional orbit consisting of sextic forms with a root of multiplicity 4 and a root of multiplicity 2 is contained in the closure of the orbit consisting of sextic forms with a root of multiplicity 4 and two distinct roots.
To show that $\brank f >2$, that is that this form $f_1=x^4y^2$ is not in the closure  of the $SL_2(\C)$-orbit of the form $f_2=x^6+y^6$, we can proceed as follow, using an idea of Bruce Reznick. To each sextic form
$f$ one can associate its catalecticant matrix $C(f)_3$ as in Example \ref{ex1}. This is an element of $Sym_4(\C)$, the vector space of symmetric $4 \times 4$ matrices with entries in $\C$. Hence we get
a well defined regular mapping $\PP(S_6) \to \PP(Sym_4(\C))$ given by
$f \mapsto C(f)_3$. The rank of the matrix $C(f)_3$ is constant when $f$ varies in an $SL_2(\C)$-orbit. Since $\rank C(f_1)_3=3$ and $\rank C(f_2)_3=2$, it follows that $f_1$ is not in the closure of the 
$SL_2(\C)$-orbit of the form $f_2$. This ends the proof of our result.
\endproof
As a consequence, we get the following analog of Corollary \ref{corGIT}, showing that the highest multiplicity $m$ of a root for a given not stable sextic $f$ determines its  border rank $\brank f$.
\begin{cor}
\label{corGITborder}
Let $f$ be a non zero binary sextic having a  root of highest multiplicity $m \geq 3$. Then the following hold.

\begin{enumerate}

\item  $m=6$ if and only if  $\brank f=1$.

\item If $m=5$, then  $\brank f=2$.

\item If $m=4$, then $\brank f=3$.

\item If $m=3$, then $\brank f=4$.

\end{enumerate}
In particular, $\rank f -\brank f =0$ for $m \in \{3,6\}$, $\rank f -\brank f =4$
for $m=5$ and $\rank f -\brank f =2$ for $m=4$.
\end{cor}

\section{Conflict of interest} 

The authors declare that they have no conflict of interest.


\begin{thebibliography}{00}



\bibitem{AH}  J. Alexander, A. Hirschowitz, Polynomial interpolation in several variables, J. Algebraic Geom. 4 (1995),  201--222.

\bibitem{BBM}  A. Bernardi, J. Brachat, and B. Mourrain, A comparison of different notions of ranks
of symmetric tensors, Linear Algebra Appl. 460 (2014), 205--230.




\bibitem{Guide} A. Bernardi, E. Carlini, M. V. Catalisano, A. Gimigliano, A. Oneto, The Hitchhiker guide to: Secant Varieties and Tensor Decomposition, Mathematics 2018, 6, 314. https://doi.org/10.3390/math6120314

\bibitem{Bo} O. Bolza, On binary sextics with linear transformations into themselves, Amer. J. Math. 10 (1888), 47--70.

\bibitem{BM} L. Brustenga i Moncus\' i, S. K. Masuti, The Waring rank of binary  binomial forms, Pacific J. Math. 313 (2021), 327--342.

\bibitem{BT} J. Buczy\' nski, Z. Teitler, Some examples of forms of high rank, Collectanea Mathematica 67(2016), 431--441.

\bibitem{B2} J. Buczy\' nski, K. Han, K., M. Mella, Z. Teitler, On the locus of points of high rank. European Journal of Mathematics 4(2018), 113--136. 

\bibitem{Car0+} E. Carlini, M.V. Catalisano, A.V. Geramita, The solution to the Waring problem for monomials and the sum of coprime monomials, J. Algebra 370 (2012) 5–14.



\bibitem{Car1+} E. Carlini, M.V. Catalisano, A. Oneto,  Waring loci and the Strassen conjecture. Adv. Math. 314(2017), 630--662.


\bibitem{Car2+} E. Carlini, M. V. Catalisano, L. Chiantini, A. V. Geramita, Y. Woo, 
Symmetric tensors: rank, Strassen’s conjecture and  e-computability,
Ann.  Scuola Normale Sup. Pisa. 18 (2018), 363--390.



\bibitem{Ci} C. Ciliberto, Geometric aspects of polynomial interpolation in more variables and of Waring’s problem. In: European Congress of Mathematics, Barcelona 2000, pages 289--316. Springer, 2001.

\bibitem{Cl}  A. Clebsch, Theorie der Bin\"aren Algebraischen Formen, Verlag von B.G. Teubner, Leipzig, (1872).


\bibitem{Co} CoCoA-5.3 (2022): a system for doing Computations in Commutative Algebra,
available at http://cocoa.dima.unige.it

\bibitem{CS} G. Comas and M. Seiguer, On the rank of a binary form. Foundations of Computational
Mathematics, 11(1)(2011), 65--78.

\bibitem{CGLM} P. Comon, G. Golub, L.H. Lim, B. Mourrain, Symmetric tensors and symmetric tensor
rank,  SIAM Journal on Matrix Analysis and Applications, 30(2008),1254--1279.



\bibitem
{Sing} { W. Decker, G.-M. Greuel, G. Pfister \and H. Sch{\"o}nemann.} \newblock {\sc Singular} {4-3-0} --- {A} computer algebra system for polynomial computations, available at {http://www.singular.uni-kl.de} (2014).



\bibitem{DStToh} A. Dimca, G. Sticlaru, Waring rank of binary forms, harmonic cross-ratio and golden ratio,  
Tohoku Math. J.,  Second Series, 74(2022).

\bibitem{Do} I. Dolgachev, Lectures on Invariant Theory, London Mathematical Society Lecture Note Series,  Cambridge: Cambridge University Press, 2003. 



\bibitem{El} E. B. Elliott, An Introduction to the Algebra of Quantics, Oxford University /
Clarendon Press, 1895. 


\bibitem{FOS} R. Fr\" oberg, G. Ottaviani,  B. Shapiro, On the Waring problem for polynomial rings Proceedings
of the National Academy of Sciences, 109 (2012), 5600--5602.

\bibitem{FLOS} R. Fr\" oberg,  S. Lundqvist, 
A. Oneto, 
B. Shapiro, 
Algebraic stories from one and from the other pockets, Arnold Math. J. 4 (2018),  137--160.



\bibitem{H} B. Hassett, Classical and minimal models of the moduli space of curves of genus two. In: Bogomolov F., Tschinkel Y. (eds) Geometric Methods in Algebra and Number Theory. Progress in Mathematics, vol 235, (2005), pp 169-192,  Birkhäuser Boston. $https://doi.org/10.1007/0-8176-4417-2_8$.



\bibitem{IK} A. Iarrobino, V. Kanev, {\it Power Sums, Gorenstein Algebras, and Determinantal Loci}, Springer Lecture Notes 1721, 1999.



\bibitem{Ig} J. Igusa, Arithmetic variety of moduli for genus two. Ann. of Math. (2) 72 1960 612–649.



\bibitem{L1} J.M. Landsberg, {\it Tensors: Geometry and Applications}, Graduate Studies in Mathematics vol.128, American Mathematical Soc., 2012.

\bibitem{LT} J.M. Landsberg, Z. Teitler, On the ranks and border ranks of symmetric tensors, Found. Comput. Math. 10(3) (2010) 339--366.






\bibitem{KSV} V. Krishnamoorthy, T. Shaska, H. V\" olklein, Invariants of Binary Forms, In: Voelklein H., Shaska T. (eds) Progress in Galois Theory. Developments in Mathematics, vol 12, pp 101--122, Springer, Boston, MA. $https://doi.org/10.1007/0-387-23534-5_6.$



\bibitem{O} A. Oneto, Waring type problems for polynomials, Doctoral Thesis in Mathematics at Stockholm University, Sweden, 2016.











\bibitem{Sy} J.J. Sylvester. Lx. on a remarkable discovery in the theory of canonical forms and of hyperdeterminants.
The London, Edinburgh, and Dublin Philosophical Magazine and Journal of
Science, 2(12)(1851), 391--410.

\bibitem{T} D. W. Taylor, Moduli of Hyperelliptic Curves and Invariants of
Binary Forms, PhD Thesis, UCLA, 2013.

\bibitem{Tok} N. Tokcan, On the Waring rank of binary forms, Linear Algebra and Its Applications 524(2017), 250--262.

	



\end{thebibliography}
\end{document}